\documentclass{conm-p-l}

\usepackage{amssymb,color}

\usepackage{graphicx}

\def\r{\mathbb R}
\def\h{\mathbb H}
\def\s{\mathbb S}
\def\c{\mathbb C}

\newtheorem{theorem}{Theorem}[section]

\newtheorem{problem}[theorem]{Problem}

\theoremstyle{definition}

\theoremstyle{remark}

\numberwithin{equation}{section}

\begin{document}

\title[Open problems on compact cmc  surfaces with boundary]{Open problems on compact constant mean curvature   surfaces with boundary}


\author{Rafael L\'opez}
\address{Departamento de Geometr\'ia y Topolog\'ia, Universidad de Granada. 18071 Granada, Spain}
\email{rcamino@ugr.es}
\thanks{}


\date{}

\begin{abstract}
We present a collection of easily stated open problems in the theory of compact constant mean curvature surfaces with boundary. We also survey the current status of answering them.

\end{abstract}

\maketitle


\section{Introduction and preliminaries}

 
In this paper, I collect a number of open problems that I have found along the years  during my research, have struggled with, and which are still open. At the same time, we present the current state of knowledge of these problems, indicating which are the difficulties in answering them.  I hope this collection of open questions may stimulate younger minds, specially because most of the problems are easily formulated and any mathematician with a basic background in classical differential geometry can understand them.

Surfaces with constant mean curvature (cmc surfaces in short or $H$-surfaces  where $H$ is the constant value of the mean curvature) belongs to the field of classical differential geometry of surfaces. A cmc surface in Euclidean $3$-dimensional space $\r^3$ has the property its area is stationary that under all volume-preserving variations in compact domains. For this reason, cmc surfaces are models of soap films and soap bubbles \cite{is} hence its interest in physics, chemistry, biology and engineering. From another perspective, the mean curvature equation is a quasilinear elliptic equation and it serves as model in the theory of PDE of elliptic type \cite{gt}.  Hence that surfaces with constant mean curvature are of interest in other fields of mathematics. As a recent sample button, overdetermined problems in PDEs initiated by Serrin, Berestycki, Cafarelli and Nirenberg \cite{bere,gnn,se2} have recently received a great boost using cmc surfaces: we refer to  \cite{pino} and references therein for the interested reader.

The main open problems     proposed in this paper comes from the three celebrated characterizations of the sphere in the family of the closed cmc surfaces. 

\begin{theorem} \label{t1}
\begin{enumerate}
\item Spheres are the only {\bf embedded} closed cmc surfaces (Alexandrov \cite{al}).
\item Spheres are the only {\bf  genus zero} cmc surfaces (Hopf \cite{ho}).
\item Spheres are the only {\bf stable} closed cmc surfaces (Barbosa-do Carmo \cite{bc}).
\end{enumerate}
\end{theorem}
There exist closed cmc surfaces others than spheres. It was Wente who firstly proved the existence of an immersed torus in $\r^3$ with constant mean curvature \cite{we3}. After that work,  new examples of closed cmc surfaces with higher genus were discovered \cite{ab,bo,ka2,ka3,ps}. 

\begin{figure}[hbtp]
\begin{center} 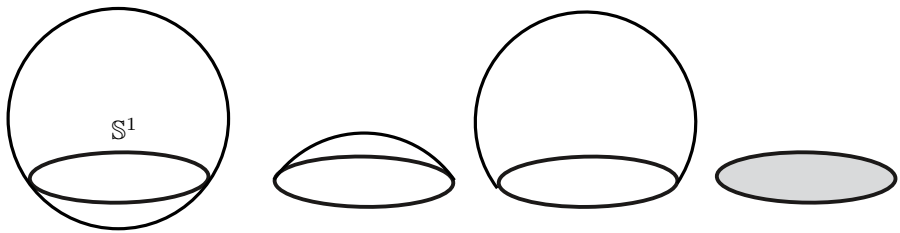\end{center}
\caption{Surfaces with constant mean curvature spanning a  circle $\s^1$ obtained by intersecting a sphere with a plane (left). The middle surfaces are (small and big) spherical caps; the right surface is a planar disk.}\label{fig1}
\end{figure}

It is natural to ask for the boundary versions of these characterizations. The simplest case of boundary curve is a circle $\s^1$.  If we intersect a sphere by a plane, we obtain two surfaces, called spherical caps, whose common boundary curve is a circle   $\s^1$   (Fig. \ref{fig1}). Both surfaces have non-zero constant mean curvature. Furthermore  the planar disk determined by $\s^1$  is a minimal ($H=0$) surface. Thus, planar disks and spherical caps are examples of cmc surfaces spanning a circle.

Following   Theorem \ref{t1}, the analogous hypothesis to consider in the boundary version is that the surface is embedded, it is a topological disk (genus zero)  or it is stable. There are the following three conjectures.  
\begin{enumerate}
\item[] \textbf{Conjecture 1.} \emph{Planar disks and spherical caps are the only compact embedded cmc surfaces in $\r^3$ with circular boundary.}
 \item[] \textbf{Conjecture 2.} \emph{Planar disks and spherical caps are the only  cmc disks immersed in $\r^3$ with circular boundary.}
\item[] \textbf{Conjecture 3.} \emph{Planar disks and spherical caps are the only compact stable cmc surfaces in $\r^3$ with circular boundary.}
\end{enumerate}
 As far as I know, the first time when similar questions appeared was in \cite{bro}: see `Problem 1.4'  proposed by  Gulliver and   Kusner. They asked about the existence of cmc surfaces spanning a circle. They proposed the hypothesis that   the surface is embedded (Conjecture 1) or it is a topological disk (Conjecture 2). For the initial question of existence of examples, Kapouleas proved that  there are compact cmc surfaces of genus $g>2$ spanning a circle \cite{ka2}. These examples appeared in the mentioned context of finding closed cmc surfaces.   

Some of the problems proposed in this paper around these conjectures  are well-known for people working in cmc surfaces, particularly, around 80 and 90's in the past century. In one way or another they appear in the works of Barbosa, Brito, Kapouleas, Koiso, Korevaar, Kusner, Meeks, Palmer, Sa Earp and Rosenberg among others.  Perhaps because their difficulty, these problems were gradually abandoned and forgotten but I think that they deserve a wider dissemination. It is also clear that the simplicity in the formulations will make famous, in this small world of differential geometry the one who can prove or disprove the above conjectures.

We use the following   notations. If $\Psi:\Sigma\to\r^3$ is an isometric immersion of a  compact surface $\Sigma$ with non-empty boundary $\partial\Sigma$ and $C\subset\r^3$ is a closed curve, we say that $\Sigma$ has $C$ as a boundary if $\Psi_{|\partial\Sigma}\colon\partial\Sigma\to C$ is a diffeomorphism. If $C$ is a planar curve, we denote by $P$ the plane containing $C$. We will assume that $P$ is horizontal, that is, parallel to the $xy$-plane when we consider $(x,y,z)$ the standard coordinates of $\r^3$. In case that $C$ is a round circle, we use $\s^1$ instead of $C$.  If the boundary  is   planar, then any minimal surface spanning $C$ is an open set of the plane containing $C$. Hence   the situation  $H=0$ can be viewed as `trivial'. In this paper we will consider surfaces in Euclidean $3$-dimensional space $\r^3$ although many of the problems can be formulated in other space forms as well as   for hypersurfaces in arbitrary dimensions.

Just like the Google reviews allow to rate the quality of a site, I have rated from 1 to 5 stars the importance, in my opinion, of each of the problems. The problems will be classified in three categories: the circular boundary case, arbitrary boundary curve and the Dirichlet problem.

\section[PP]{Part I:   circular boundary}\label{sec2}

In this section we will investigate the case when the boundary of the surface is a circle $\s^1$.  Spheres are rotational surfaces.  Rotational surfaces of constant mean curvature are called Delaunay surfaces and they are well known: planes, catenoids, spheres, cylinders, unduloids and nodoids \cite{de}. The profile curves of the  Delaunay surfaces are the trace of the focus of a conic which is rolled without slipping along a line and rotating this curve around that line.  Planes and catenoids are the only ones with $H=0$. Taking compact pieces of Delaunay surfaces, planar disks and  spherical caps are the only compact rotational cmc surfaces spanning a circle.

The first problem is the boundary version of the Alexandrov theorem. 
\begin{problem}[5*]\label{p1}
Is a compact embedded cmc surface of $\r^3$ spanning a circle a planar disk or a spherical cap?
\end{problem} 

Let us show the difference from the closed case. The technique employed in the Alexandrov theorem is called the reflection method \cite{al}. If $\Sigma\subset\r^3$ is a closed embedded surface, then   $\Sigma$ encloses a $3$-domain $W\subset \r^3$. Let $N$ be the orientation pointing towards $W$. Fix an arbitrary direction of $\r^3$, say the $x$-direction. Let $\Pi_t$ be the plane of equation $x=t$. Denote $\Sigma_t^+=\Sigma\cap\{x>t\}$, $\Sigma_t^{-}=\Sigma\cap\{x<t\}$ and $\Sigma_t^*$ the reflection of $\Sigma_t^{+}$ with respect to $\Pi_t$. Coming from $t=+\infty$ towards $\Sigma$, let $\Pi_{t_0}$ be the first plane touching $\Sigma$. By embeddedness of $\Sigma$,  $\Sigma_t^{*}\subset W$ for all $t\in (t_0-\epsilon,t_0)$ and for some $\epsilon>0$. Decreasing $t$ and by compactness, let $t_1<t_0$ be the first time where $\Sigma_{t_1}^*$ touches $\Sigma_{t_1}^-$. If $p\in\Sigma_{t_1}^{*}\cap \Sigma_{t_1}^-$, then the surfaces $\Sigma_{t_1}^*$ and $\Sigma_{t_1}^-$ are $H$-surfaces with the same orientations  at $p$. Moreover $\Sigma_{t_1}^*$ lies in one side of $\Sigma_{t_1}^-$ around $p$. Then the  maximum principle  implies $\Sigma_{t_1}^*=\Sigma_{t_1}^-$ in an open set around $p$ \cite{gt}. By connectedness, $\Sigma_{t_1}^*=\Sigma_{t_1}^-$. This shows that $\Pi_{t_1}$ is a plane of symmetry of $\Sigma$. Since the fixed orientation was arbitrary, then $\Sigma$ is a sphere.
 
 Consider now a compact embedded cmc surface spanning a circle $\s^1$. Suppose that $\s^1\subset P$ is  centered at the origin. If we repeat the reflection  argument, a first problem   we face is  the existence of the domain $W$.  For example, if we simply attach to $\Sigma$ the planar disk $D$ bounded by $\s^1$, the existence of $W$ is not assured because the surface may across $D$.    

\begin{figure}[hbtp]
\begin{center} 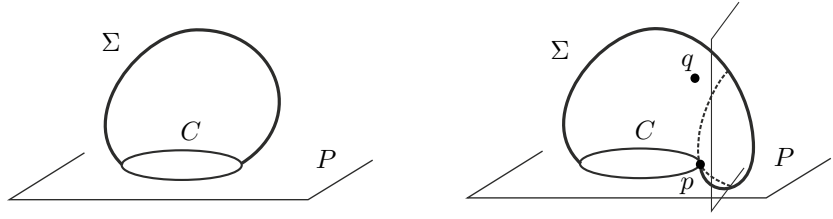\end{center}
\caption{Left: a compact embedded surface over the boundary plane $P$. Right: an embedded  surface where the Alexandrov reflection method fails at $p$. Notice that the surface $\Sigma$ is star-shaped from the point $q$.}\label{fig2}
\end{figure}

A situation where    the reflection method applies is when $\Sigma$ lies in one side of $P$.  This hypothesis is not only because now $\Sigma\cup D$ encloses a $3$-domain $W$, but also because $\Sigma$ does not intersect the plane $P$: see Fig. \ref{fig2} left. By repeating the argument as in the closed case, we now fix an arbitrary horizontal direction, say the $x$-direction. The time $t_1$ can occur when $p\not\in\s^1$ or $p\in\s^1$. In the first case, then $\Pi_{t_1}$ is a plane of symmetry, in particular, of symmetry of $\s^1$. This implies $t_1=0$. The second case implies $t_1=0$. Let us follow  this situation. If $\Sigma_0^*$ is tangent to $\Sigma_0^{-}$, then the  maximum principle applies proving that $\Pi_0$ is a plane of symmetry.  On the contrary, $\Sigma_0^*$ is contained in $W$ and it does not touch $\Sigma_0^-$. We see that this situation is not possible. For this we repeat the argument but coming from $t=-\infty$. Then there is a time $t_2<0$ such that the reflection of $\Sigma_{t_2}^{+}$ about $\Pi_{t_2}$ touches tangentially $\Sigma_{t_2}^{-}$. By the maximum principle, the plane $\Pi_{t_2}$ would be a plane of symmetry of $\Sigma$. This is a contradiction because $t_2<0$ and $\Pi_{t_2}$ cannot be a plane of symmetry of $\s^1$. Summarizing, we have proved that for any horizontal direction, the plane orthogonal to this direction and containing the $z$-axis is a plane of symmetry of $\Sigma$. This proves that $\Sigma$ is rotationally symmetric about the $z$-axis, hence a spherical cap.

Thus the Problem \ref{p1} can be affirmatively answered if one proves that   a compact embedded cmc surface spanning  a circle is contained in one side of $P$. 

\begin{problem}[3*]\label{p2}
Under what conditions on the boundary,  a compact embedded cmc surface of $\r^3$ spanning a planar closed curve $C$ is contained  in one side of the boundary plane? and  if $C$ is convex?
\end{problem} 

A partial answer to Problem \ref{p2} was obtained by Koiso \cite{ko} assuming that $\Sigma\cap (P\setminus\overline{D})=\emptyset$. In fact,   in \cite{ko} it is not necessary that $H$ is constant but only $H\not=0$ on $\Sigma$.

  The existence of the domain $W$ under other circumstances does not solve Problem \ref{p1}. In Fig. \ref{fig2} right, the surface $\Sigma$ together with $D$ defines $W$ but $\mbox{int}(\Sigma)$ intersects the plane $P$. Then  the reflection process may find  a touching point $p$ between   $\Sigma_t^*$ and $\partial\Sigma$, where $\Sigma_t^*$ and $\Sigma_t^-$  are not tangent: the maximum principle cannot be applied. 

Let us observe that  the surface   of Fig. \ref{fig2} right  is star-shaped with respect to the point $q$, that is,   any ray from $q$ meets the surface once at most. If $\Sigma$ is compact, the star-shaped property implies that $\Sigma$ is embedded. This is equivalent to say that $\Sigma$ is a radial graph on some domain of a sphere centered at $q$.    In the XIXth century,  Jellet proved that a closed star-shaped cmc surface is a sphere \cite{je}. This motivates the boundary version of this result. 

\begin{problem}[5*]\label{p3}
Is a compact star-shaped cmc surface of $\r^3$ spanning a circle a planar disk or a spherical cap?
\end{problem} 
Although a positive answer   is a partial answer to Problem \ref{p1}, I rate it as a 5* problem because it is the analog of a classical result. The proof in the closed case is the following. Assume that $\Sigma$ is star-shaped from the origin of $\r^3$. If $x$ is the position vector of $\Sigma$,  the support function $h=\langle N,x\rangle$ has sign, say $h<0$ with the inward orientation (and thus, $H>0$). Moreover we have
\begin{equation}\label{e1}
\begin{split}
\Delta |x|^2&=4+4H h\\
 \Delta h&=-2H-|\sigma|^2h,
\end{split}
\end{equation}
where $\Delta$ is the Beltrami-Laplace operator and $\sigma$ is the second fundamental form of $\Sigma$. The first equation holds for any surface but   the second one holds because $H$ is constant.   Integrating and using that $\Sigma$ is closed, we have 
$A+H\int_\Sigma h=0$ and $2HA+\int_\Sigma|\sigma|^2h=0$, where $A=\mbox{area}(\Sigma)$. Combining both equations and using $|\sigma|^2\geq 2H^2$, we have
$$0=2HA+\int_\Sigma |\sigma|^2h\leq 2HA+2H^2\int_\Sigma h=2H(A+H\int_\Sigma h)=0.$$
Thus $|\sigma|^2=2H^2$ on $\Sigma$, concluding that the surface is totally umbilical, hence a sphere. 
 In the boundary case, the troubles appear in  the boundary terms coming from the divergence theorem for both equations \eqref{e1}.

Returning to the difficulty to find the domain $W$, it was proved in \cite{bemr} that if   $\partial\Sigma=C$ is   convex   and   $\Sigma$  is transverse to the boundary plane $P$ along $C$, then $\Sigma$ lies in one side of $P$. Thus in Problem \ref{p1} it would remain to study that it is not possible that along $\s^1$, the surface have parts in both sides of $P$.  In the same paper, the authors propose the following  problem. 

\begin{problem}[4*]\label{p4}   Let $C$ be a simple closed convex curve in a plane $P$. Let $\Sigma$ be a compact embedded cmc surface in $\r^3$ with $\partial \Sigma=C$. Assuming that $\Sigma$ lies in one side of $P$, is $\Sigma$   a topological disk?
\end{problem} 

 A partial  affirmative answer is when  $\Sigma$ is a graph on $D$ in a neighborhood of $C$. Under this hypothesis,  $\Sigma$ is a graph on $D$.  Here $D\subset P$ is the domain bounded by $C$ \cite{lo2}. In this result is not necessary that $C$ is convex.  The proof is an application of the Alexandrov method by means of parallel planes to $P$. Returning to  Problem \ref{p4} we propose:  
 
 \begin{problem}[4*]\label{p5}   Under the same hypothesis of Problem \ref{p4}, if $\Sigma$ is a graph on $P\setminus\overline{D}$ around $C$, is $\Sigma$ a bigraph? that is,   is there a planar domain $\Omega$, $\overline{D}\subset\Omega$, such that   $\Sigma=\Sigma_1\cup \Sigma_2$, where $\Sigma_1$ is a graph on $D$ and $\Sigma_2$ is an annulus which it is a graph on  $\Omega\setminus\overline{D}$?
 \end{problem}

Related with the above two problems, Ros and Rosenberg answered affirmatively  Problem \ref{p4} when the volume of $\Sigma$ is sufficiently large\cite{rr}.

A further related question     is when the  boundary $\partial\Sigma$ is formed of two curves. 

 \begin{problem}[4*]\label{p6}    Let $C_1$ and $C_2$ be two convex curves in parallel planes. Is a   compact embedded cmc surface  spanning $C_1\cup C_2$   a topological annulus?
 \end{problem}

If $H=0$, this problem was proposed by Meeks \cite{me} and it is still open. Problem \ref{p6} has an interesting version when the boundary is formed by two coaxial circles.  

\begin{problem}[4*]\label{p7} Let $C_1$ and $C_2$ be two coaxial circles in parallel planes. Is a   compact embedded cmc surface  spanning $C_1\cup C_2$   a surface of revolution?
\end{problem}

 In order to apply the Alexandrov reflection method,  if the surface is contained in the slab  determined by the two boundary planes, then the surface is rotational \cite{lo1}. However, and in contrast to  Problem \ref{p2}, one cannot expect that the surface is included in the slab: there are embedded pieces of nodoids spanning two coaxial circles with parts outside of the slab.  
  
We now continue with Conjecture 2.

\begin{problem}[5*] \label{p8}
Is a topological cmc disk spanning a circle a planar disk or a spherical cap?
\end{problem}

In the closed case, the proof is based in the properties of the Hopf differential form.  Assume that $\Psi=\Psi(z)$ is a conformal parametrization of the surface, where $z=u+iv$. Now we denote by $\langle,\rangle$   the metric in $\r^3$ as well as in $\c^3$. Then  $2\langle\Psi_z,\Psi_{\bar{z}}\rangle=E>0$ and $4HE=\langle N,\Psi_{uu}+\Psi_{vv}\rangle$. Hopf defined the quadratic form   $Q(dz)^2$ where $Q=\langle N,\Psi_{zz}\rangle=\frac14(e-g+2if)$, being $\{e,f,g\}$ the coefficients of $\sigma$ with respect to $\Psi$. Note that the zeroes of $Q$ are the umbilical points of $\Sigma$. Using the Codazzi equations, one gets $Q_{\bar{z}}=EH_z/2$ \cite{ho}. Consequently  $Q$ is holomorphic if and only if $H$ is constant. Moreover either the   umbilical points of the surface are isolated or the surface is totally umbilical. Suppose  now that $\Psi$ is an immersion of a topological sphere. Then we can assume that the immersion is $\Psi\colon\overline{\mathbb{C}}\to\r^3$. Because $Q\colon\overline{\mathbb{C}}\to\mathbb{C}$ is holomorphic, then $Q\equiv 0$. Thus $\Psi$ is totally umbilical, hence it is a sphere. In the boundary version, if we   have a cmc topological disk $\Psi$, the zeroes of $Q$ coincide again with the umbilical points of $\Psi$. However,  we cannot assure   $\Psi\equiv 0$ on the surface because the existence of boundary of $\Sigma$.

A result related to Problem \ref{p8} is the following  \cite{lm1}. For immersed disks $\Sigma$ in $\r^3$    Barbosa and do Carmo proved the isoperimetric inequality  $L^2\geq 4\pi A(1-\frac{k_0A}{4\pi})$, where  $L$ is the length of $\partial \Sigma$, $A$ is the area of $\Sigma$ and $k_0\in\r$ is an upper bound of $K$ \cite{bc1}. The equality holds if and only if $\Sigma$ is a totally geodesic disk with curvature $K\equiv k_0$. If $\Sigma$ has constant mean curvature $H$, we can choose $k_0=H^2$. In case that $\partial\Sigma$ is a circle $\s^1(r)$ of radius $r$,   this inequality implies that either $A\leq A_{-}$ or $A\geq A_{+}$, where  
$$A_{-}= \frac{2\pi}{H^2}(1-\sqrt{1-r^2H^2)},\quad A_{+}= \frac{2\pi}{H^2}(1-\sqrt{1+r^2H^2)}.$$
The values $A_{-}$ and $A_{+}$ are the areas of the small and the big $H$-spherical caps with boundary $\s^1(r)$, respectively (Fig. \ref{fig1}). In \cite{lm1} it is proved that if $\Sigma$ is a topological $H$-disk with $A\leq A_{-}$, then $\Sigma$ is the small spherical cap.

The following problem is the Conjecture 3.

\begin{problem}
 Is a stable compact cmc surface spanning a circle a circle a planar disk or a spherical cap?
\end{problem}

Stability means that the second derivative of the area functional is non-negative for any  volume preserving variation of the surface. This is equivalent to 
\begin{equation}\label{es}
-\int_\Sigma u(\Delta u+|\sigma|^2 u)\, \geq 0
\end{equation}
for all $u\in C^2(\Sigma)$ such that 
\begin{equation}\label{es2}
\int_\Sigma u\,  =0,\quad u_{|\partial\Sigma}=0.
\end{equation}
The mean zero integral $\int_\Sigma u\,  =0$ is a consequence of the fact that all admissible  variations of $\Sigma$  preserve the volume of the surface. The boundary condition $u_{|\partial\Sigma}=0$ says that the variation fixes the boundary of the surface. In the closed case the condition $u_{|\partial\Sigma}=0$   makes no sense and it is dropped. In such a case,   Barbosa and do Carmo consider the function $u=1+H\langle N,x\rangle$   in \eqref{es} as from \eqref{e1},   we know that $\int_\Sigma u=0$. After some computations, and using the second identity of \eqref{e1},   the stability  inquality \eqref{es} becomes
$$-\int_\Sigma (|\sigma|^2-2H^2)\, \geq 0.$$
Thus   $|\sigma|^2=2H^2$, hence $\Sigma$ is totally umbilical, and thus, it is a sphere.  

Coming back to the boundary case, the extra condition $u_{|\partial\Sigma}=0$ adds difficulties to find a suitable test function to put in \eqref{es}. Notice that the above function $u=1+H\langle N,x\rangle$ does not satisfy  $\int_\Sigma u\,  =0$ neither $ u_{|\partial\Sigma}=0$. In \cite{alpm} it was proved that if $\Sigma$ is a stable topological cmc disk, then it is a planar disk or a spherical cap. This gives a partial answer to Conjecture 3. The proof uses the Hopf differential $2$-form $Q(dz)^2$ and the behavior of the first and second Dirichlet eigenfunctions of the Jacobi equation.

To finish this section, we return to  Kapouleas' examples \cite{ka2}.   These examples have self-intersections, genus $g>2$ and they are included in one of the half-spaces determined by the boundary plane. The methods employed in \cite{ka2} are the same that for  the construction of closed cmc surfaces. There are a great number of  pictures of closed cmc  surfaces: see for example, \cite{g1,g2,g3}. However, up today, we do not have any numerical pictures of the Kapouleas' examples with circular boundary.

\begin{problem}[3*]  Give numerical pictures of the Kapouleas' examples with circular boundary.
\end{problem}

As far as we know,  graphics of cmc surfaces with circular boundary a circle appear in \cite{brd}, where the authors solve the Bj\"{o}rling problem for $H$-surfaces.  In  Fig. 1 of \cite{brd}, the authors take a circle $\s^1$ in the Bj\"{o}rling data and the graphic shows the corresponding Bj\"{o}rling solution. Of course, the result is local in the sense that the surface is defined in a neighborhood of $\s^1$: it is unknown if the surface closes to be compact.

Related with the three conjectures, an easier  approach is assuming  two of the three hypothesis. Recall that only if the surface is a topological disk and   stable, it is known that the surface is a planar disk or a spherical cap \cite{alpm}. 

\begin{problem}[4*] Is a compact  stable embedded cmc surface of $\r^3$ spanning a circle a planar disk or a spherical cap?
\end{problem}

\begin{problem}[4*] Is a compact    embedded cmc disk of $\r^3$ spanning a circle a planar disk or a spherical cap?
\end{problem}

  The   proposed problems in this section can be formulated in space forms.  Notice that examples of compact cmc surfaces spanning a circle are geodesic spheres of totally umbilical surfaces.  The   hyperbolic space $\h^3$  is special  because the behavior of the $H$-surfaces in $\h^3$ depends on the value of $H$ in relation with the value $1$, the absolute value of the sectional curvature of $\h^3$. If $|H|=1$, then $1$-surfaces are similar to minimal surfaces of $\r^3$. If $|H|<1$, then $H$-surfaces have not a counterpart in $\r^3$. In this case, the three conjectures are true. Even more, any compact $H$-surface in $\h^3$ spanning a circle is umbilical without any further assumption on the surface \cite{be2,lo7}. However, if $|H|>1$, then the properties of $H$-surfaces are similar to cmc surfaces of $\r^3$:  the same results that hold in $\r^3$ are valid in $\h^3$; the same open questions in $\r^3$ are open in $\h^3$.

  In Lorentzian space forms, the mean curvature is  defined for surfaces with non-degenerate induced metric. However the behavior of spacelike surfaces and timelike surfaces is very different. If  the ambient space is the Lorentz-Minkowski space, it was proved that the only immersed compact spacelike cmc surfaces spanning a circle are   planar disks and   hyperbolic caps: no further assumptions are needle \cite{alpastor}.

\section{Part II:  arbitrary boundary curve}

In this section we will consider that the boundary   of a compact cmc surface is not necessarily a circle.  
 Let $C\subset\r^3$ be a  closed curve. If $\Sigma$ is a compact $H$-surface bounded by $C$, then the value of $H$ is not arbitrary and it depends on the geometry of $C$. This is due to the following argument (\cite{lm2}; also in     \cite{he2,kk,ku} for particular cases). The $1$-form $v\mapsto (H   p+N)\times v$, $v\in T_p\Sigma$, is closed because $H$ is constant. The divergence theorem yields $H\int_{\partial\Sigma}\alpha\times\alpha'+\int_{\partial\Sigma}\nu=0$, where $\alpha$ is a parametrization by arc-length of $\partial\Sigma$ and $\nu=N\times\alpha'$ is the unit conormal vector along $\partial\Sigma$. If $\vec{v}\in\r^3$ is a unit vector, then 
 \begin{equation}\label{flux}
 H\int_{\partial\Sigma} \langle \alpha\times\alpha',\vec{v}\rangle+\int_{\partial\Sigma}\langle \nu,\vec{v}\rangle=0.
 \end{equation}
 This identity is know as the flux formula of cmc surfaces because $\int_{\partial\Sigma}\langle \nu,\vec{v}\rangle$ is the flux of the vector field $\vec{v}$ along $\Sigma$. For example, on a compact minimal surface, the the sum of the flux along all the boundary components must vanish.
 
 The first integral   $\int_{\partial\Sigma} \langle \alpha\times\alpha',\vec{v}\rangle$ is twice the algebraic area $ {A}_{\vec{v}}$ of the projection of $C$ on    a plane $\Pi$ orthogonal to $\vec{v}$. Then $2|H||{A}_{\vec{v}}|\leq\int_{\partial\Sigma}|\langle\nu,\vec{v}\rangle|\leq L$, where $L$ is the length of $C$. If $H\not=0$, we conclude 
 $$|H|\leq \frac{L}{2|{A}_{\vec{v}}|}.$$
 This implies that $H$ is not arbitrary and that the value of $H$  only depends on $C$.    For example, suppose that  $C$ is a simple closed curve contained in a plane $P$ and let $D$ be the bounded domain by $C$ in $P$. Then $\max\{|{A}_{\vec{v}}|:|\vec{v}|=1\}$ is attained when $\vec{v}$ is orthogonal to $P$ and its value is the area of $D$. Therefore    $|H|\leq \frac{L}{2|D|}$.   If $C=\s^1(r)$, then  $|H|\leq\frac{1}{r}$.  The case $\s^1(r)$ as boundary curve is special because for each   $-\frac{1}{r}\leq H\leq\frac{1}{r}$ there exists an $H$-surface spanning $\s^1(r)$: planar disks and spherical caps. This phenomenon is unknown when $C$ is an arbitrary curve. 
\begin{problem}[3*]\label{pr31} Let $C$ be a  simple planar closed curve.  Is there an $H$-surface bounded by $C$ for each   $H$ of the interval $[0,\frac{L}{2|D|}]$? Particularize to the case that $C$ is a convex curve or   $C$ is an ellipse.
\end{problem}

If $J_C$ is the subset of $\r$ formed by the numbers $H$ such that there is an $H$-surface bounded by $C$,   Problem \ref{pr31} asks the relation between $J_C$ and the interval $[0,\frac{L}{2|D|}]$, where we know that $J_C\subset [0,\frac{L}{2|D|}]$. 
It is also known that for values close to $0$, there exist $H$-graphs on $D$ bounded by $C$ (independently whether $C$ is convex or not). So, there is $\epsilon>0$ such that $[0,\epsilon)\subset J_C$. Nothing is   known about the distribution of values $H$ in the set $J_C$. The last value $\frac{L}{2|D|}$ is of interest. If $C=\s^1(r)$, this value is $1/r$. It was proved that hemispheres are the only $1/r$-surfaces bounded by $\s^1(r)$: not further hypotheses are needed \cite{be}. For a simple planar closed curve, if $H=\frac{L}{2|D|}$, the flux formula \eqref{flux} says that $\Sigma$ is orthogonal to $P$ along $C$. In particular, $C$ is a line of curvature of $\Sigma$ by the Joachimsthal theorem. With the extra hypothesis that  $C$ is convex and $\Sigma$ is embedded, we know from  \cite{bemr} that $\Sigma$ lies in one side of $P$. Then the Alexandrov reflection method can be employed  obtaining that $C$ is a circle and $\Sigma$ is a spherical cap. This result also holds assuming that the angle between $\Sigma$ and $P$ is constant (but not necessarily $\pi/2$). In the general case, if $\Sigma$ is an immersed $H$-disk and $C$ is a simple closed curve, then $|H|=\frac{L}{2|D|}$ implies that   $\Sigma$ is   a hemisphere and $C$ is a circle \cite{lo10}. 

\begin{problem}[3*] Let $C$ be a simple planar closed curve. Is there a $H$-surface bounded by $C$ where $H=\frac{L}{2|D|}$? What can be said about the geometry of this surface?
\end{problem}

A method to construct    cmc surfaces is by solving the Plateau problem: given a boundary curve $C$ and $H\in\r$, find an immersed disk of mean curvature $H$ and boundary $C$. The techniques of functional analysis look for a critical point of the functional $X\mapsto \int_{\Omega}\frac12|\nabla X|^2+\frac{2H}{3}\int_{\Omega}\langle X,X_u\times X_v\rangle$, where $X=X(u,v) \colon \Omega\to\r^3$ is a conformal map from the unit disk $\Omega$  and $X_{|\partial \Omega}$ a parametrization of $C$. Critical points of this functional are $H$-surfaces.   Hildebrandt proved  that  if $C$ lies included in a ball $\mathcal{B}_R$ of  radius $R>0$, then for any $H$ such that  $|H|<1/R$, there is an $H$-surface bounded by $C$. Moreover, this surface is included in the closure of $\mathcal{B}_R$ \cite{hi1}. Later, Gulliver and Spruck proved a similar result replacing $\mathcal{B}_R$ by a solid cylinder $\mathcal{C}_R$ of radius $R$ \cite{gs1}. It is natural to ask if it is possible to get a similar result bounding only  one spatial coordinate. Notice that given a circle $\s^1(r)$, $H$-spherical caps bounded by $\s^1(r)$ lies in a slab of width $2/|H|$.

\begin{problem}[3*] Let  $H\not=0$ and let $C$ be a closed curve contained in    a slab $\mathcal{S}_w\subset\r^3$ of width $w$ such that  $|H|<2/w$. Does  an $H$-topological disk exist spanning $C$   and included in $\mathcal{S}_w$?
 \end{problem}
 
 It is crucial in the results \cite{gs1,hi1}  that the boundary of the sphere or of the cylinder has positive mean curvature: this is not the case of a slab.  Related with  the results of  Hildebrandt, Gulliver and Spruck,   Barbosa considered the particular case that $C$ is a circle. He proved that   an $H$-surface bounded by a circle $\s^1(r)$ and included in a solid cylinder of radius $1/|H|$ must be a spherical cap \cite{ba1,ba2}: here there is no further hypothesis on the surface. Thus we have the following problem
 
 \begin{problem}[3*] Let $H\not=0$ and let $\Sigma$ be  a compact $H$-surface spanning   a circle $\s^1(r)$. If $\Sigma$ included in the slab $\mathcal{S}_w$, with $|H|<\frac{2}{w}$, is $\Sigma$ a spherical cap? If we replace  $ \frac{2}{w}$ by $w\frac{1}{w}$, is $\Sigma$ a small spherical cap?  
 \end{problem}

The results proved in \cite{ba1,ba2} also hold for closed cmc surfaces. However, in the case that the surface is included in a slab, we can propose:
\begin{problem}[3*]   
Is the sphere of radius $1/H$  the only closed $H$-surface    included in the closure of a slab of width $2/|H|$?
\end{problem}

Notice that small spherical caps bounded by $\s^1(r)$ are included in the solid cylinder, where $D\subset P$ is the domain bounded by $\s^1(r)$.

\begin{problem}[3*]   Let $C\subset P$ be a simple planar convex curve.  Let $\Sigma$ be a compact cmc surface spanning $C$.  If  $\Sigma$ is included in  $D\times\r$, is $\Sigma$   a graph on $D$?
\end{problem}

If $\Sigma$ is embedded, then $\Sigma$ is a graph \cite{lo2}. On the other hand, if there exists an $H$-graph $\Sigma^*$ on $D$, where $H$ is the mean curvature of $\Sigma$, then $\Sigma=\Sigma^*$ up to reflections about $P$ \cite{lo10}. This is a consequence of the flux formula \eqref{flux} and the maximum principle. In particular, this is the situation when $H$ is sufficiently small because  the existence of $H$-graphs on $D$ with boundary $C$ is assured when $H$ varies in an interval $[0,\epsilon)$ for some $\epsilon>0$.  The above   two results does not  use the convexity of $C$.

In the same direction of the   results in \cite{gs1,hi1},  Wente proved that given   a real number $V>0$, there is a topological  disk $\Sigma$ bounded by $C$ and with volume $V$ \cite{we2}. This  surface is a  critical point of the Dirichlet energy $ \int_\Omega\frac12|\nabla X|^2$   in a suitable family of conformal immersions $X\colon \Omega\to\r^3$ with boundary   $C$ and fixed volume $V$. In particular, the mean curvature of $\Sigma$ must be constant.  In case that $C$ is a   convex planar curve, it was proved in \cite{lm1} that a cmc surface bounded by $C$ is a graph if the volume of $\Sigma$ is sufficiently small. Recall here the   result  of Ros and Rosenberg for large volume  \cite{rr}, mentioned in the second section.

The following problem was first formulated by Rosenberg  \cite[p. 21]{hoff}.

\begin{problem}[4*] \label{pr35}Let $C$ be a non-circular closed  convex planar curve. Prove that there is   $V_C>0$ depending only on $C$ such that, if $\Sigma$ is a compact cmc embedded disk   spanning $C$ and with volume $V>0$, it holds:
\begin{enumerate}
\item If $V\leq V_C$, then $\Sigma$ lies in one side of  $P$.
\item If $V=V_C$, then $\Sigma$ is tangent to $P$ at the points of minimum curvature of $C$.
\item If $V>V_C$, then $\Sigma$ has points in both sides of $P$.
\end{enumerate}
\end{problem}

Hoffman used the Brakke's Surface Evolver to show numerical graphics that support that (3) is true  \cite{hoff}. Notice that   $C$ is not a circle in  Problem \ref{pr35}.

When the volume is sufficiently small and motivated by the result of \cite{lm1}, we can ask the following problem.  

\begin{problem}[4*] 
  Let $C_1,C_2\subset\r^2$ be two simple planar convex curves with disjoint interiors.   Does   $V_0>0$ exist such that if  $\Sigma$ is a cmc surface with boundary  $C_1\cup C_2$ and volume $V<V_0$, then $\Sigma$ is not connected.
\end{problem}

\section{The Dirichlet problem for the cmc equation}\label{sec4}

A method  to find $H$-surfaces spanning a given curve is by solving the  Dirichlet problem
\begin{equation}\label{eh}
\left\{
\begin{split}
\mbox{div}\frac{Du}{\sqrt{1+|Du|^2}}&=2H,\quad \mbox{in $\Omega$}\\
u&=\varphi,\quad\mbox{in $\partial\Omega$},
\end{split}\right.
\end{equation}
where  $  u\in C^2(\Omega)\cap C^0(\overline{\Omega})$. Here  $\Omega$ is a domain of $\r^2$, $H\in\r$ and  $\varphi:\partial\Omega\to\r$ is a continuous function.   The graph $\Sigma=\{(x,y,u(x,y)):(x,y)\in\Omega\}$ is an $H$-surface whose boundary $C$ is the graph of $\varphi$. Notice that integrating the first equation of \eqref{eh} and using the divergence theorem, we obtain the flux formula \eqref{flux}. By the maximum principle, the solution of \eqref{eh}, if exists, is unique. The solvability of \eqref{eh} requires to establish  {\it a priori} $C^1$ estimates for the solutions of \eqref{eh}. This means finding     estimates  $|u|$  in $\Omega$ and gradient estimates  $|Du|$   along $\partial \Omega$. In general,    results in the literature employ   suitable $H$-surfaces   as barriers to get these estimates.   A result of Serrin shows that \eqref{es} has a solution for {\it any } function $\varphi$ if and only if $\Omega$ is convex and the curvature $\kappa$ of $\partial\Omega$ satisfies $\kappa>2|H|$ \cite{se1}. 
 
 From now on we assume $\varphi=0$ (planar boundary curve). Then  it is natural to relax the Serrin's condition on $H$. For example, if $\Omega$ is a convex domain, a first result asserts that if   $\kappa>|H|$ then   \eqref{eh} is solvable. Here the condition  $\kappa>|H|$ allows the use of suitable  small $H$-spherical caps as barriers in establishing the $C^1$-estimates \cite{lo10}.  In \cite{lm2} it was proved an  estimate of the height   of an $H$-surface $\Sigma$ with planar boundary, namely,
 \begin{equation}\label{area} 
 h\leq \frac{2\pi}{|H|}A,
 \end{equation}  
 where $h$ is the height of $\Sigma$ with respect to the boundary plane and  $A$ is the area of $\Sigma$. Thanks to \eqref{area}, the Dirichlet problem   \eqref{eh} was solved under  a variety of assumptions, namely,  $2|\Omega|H^2<\pi$, $LH^2<\sqrt{3}\pi$ and $\mbox{diam}(\Omega)<1/|H|$ \cite{lo10,lo13,lo14,lm1}; see also \cite{lo15,ri,ri2} for other results of existence. However, if $\Omega$ is a round disk, the above inequalities are not sharp such as small spherical caps show. 

\begin{problem}[4*] \label{p41}Let $\Omega$ be a bounded convex domain. Solve the Dirichlet problem \eqref{eh} for $\varphi=0$ under each one of the following hypothesis:
\begin{enumerate}
\item $|\Omega| H^2<\pi$.
\item $LH^2<2\pi$.
\item $\mbox{diam}(\Omega)<2/|H|$.
\end{enumerate}
\end{problem}
If $\Omega$ is a non-bounded convex domain, it was proved in \cite{lo13} that  \eqref{eh} has a solution for $\varphi=0$  if and only if $\Omega$ is included in a strip of width $1/|H|$. Notice that a half-cylinder of radius $\frac{1}{2|H|}$ is a graph on a strip of width $1/|H|$. Interesting questions on cmc graphs   on convex domains are about the convexity of the solutions or of its level curves. Here we are assuming planar boundary  curve ($\varphi=0$ on $\partial\Omega$). It is known that the level curves are convex if the mean curvature $H$ or the volume  $V$ of the graph are sufficiently small \cite{mcc,sa}. However,  Wang gave an example of a cmc graph on a convex domain whose level curves are not all convex \cite{wan}.

\begin{problem}[3*] Suppose that $\Omega$ is a bounded convex domain. Let $u$ be a solution \eqref{eh} for $\varphi=0$ on $\partial\Omega$. Give a control of the distribution of the convex level curves in terms of the value $H$. Extend the Wang's example, if possible, when $\Omega$ is not bounded.
\end{problem}

Following   \eqref{eh} and $\varphi=0$ along $\partial\Omega$, it is natural to ask if a solution of the Dirichlet problem is convex in the sense that the Gauss curvature $K$ of the graph is positive everywhere. Notice that there exist pieces of unduloids obtained by cutting-off an unduloid for some particular planes parallel to the rotation axis that are graphs on convex domains but the surface has parts with negative Gauss curvature.   
See related results in  \cite{den,hua,lo-gra}.

 \begin{problem}[4*] Let $\Omega$ be a bounded convex domain. Give conditions that assure the convexity of  a cmc graph on $\Omega$   and with planar boundary \end{problem}

\section*{Acknowledgments} 
The author is  a member of the IMAG and of the Research Group ``Problemas variacionales en geometr\'{\i}a'',  Junta de Andaluc\'{\i}a (FQM 325). This research has been partially supported by MINECO/MICINN/FEDER grant no. PID2023-150727NB-I00,  and by the ``Maria de Maeztu'' Excellence Unit IMAG, reference CEX2020-001105- M, funded by MCINN/AEI/10.13039/501100011033/ CEX2020-001105-M. Finally, I thank the editors of the CONM volume  for the possibility to contribute.


 \end{document}